\documentclass[12pt,reqno]{amsart}

\usepackage{amsmath,amssymb,amsthm}
\usepackage{graphicx}
\usepackage[all]{xypic}
\usepackage[pdftex]{hyperref}
\input xypic
\newtheorem{problem}{Problem}[section]
\newtheorem{example}{Example}[section]
\newtheorem{remark}{Remark}[section]
\newtheorem{theorem}{Theorem}[section]
\newtheorem{corollary}{Corollary}[section]
\newtheorem{lemma}{Lemma}[section]                                                                                                                                                                                      \newtheorem{Proposition}{Proposition}[section]
\newtheorem{definition}{Definition}[section]

\numberwithin{equation}{section}

\def\rem{\begin{remark}}
\def\erem{\end{remark}}
\def\ex{\begin{example}
  }
  \def\eex{\end{example}}
\def\thr{\begin{theorem}}
\def\ethr{\end{theorem}}
\def\pro{\begin{Proposition}}
\def\epro{\end{Proposition}}
\def\coro{\begin{corollary}}
\def\ecoro{\end{corollary}}
\def\df{\begin{definition}}
\def\edf{\end{definition}}
\def\lm{\begin{lemma}}
\def\elm{\end{lemma}}
\def\pf{\begin{proof}}
\def\epf{\end{proof}}
\def\problem{\begin{problem}}
\def\eproblem{\end{problem}}

\def\it{\begin{itemize}}
\def\hit{\end{itemize}}
\def\equ{\begin{equation}}
\def\eequ{\end{equation}}

\newcommand{\vm}{\forall}

\newcommand{\mtn}{\rightarrow}

\def\Hom{\mbox{Hom}}

\def\r{\Bbb R}

\def\z{\Bbb Z}

\begin{document}
\title[Jensen's functional equation on the symmetric group $\bold{S_n}$ ]
{Jensen's functional equation on the symmetric group $\bold{S_n}$}

\author{C\^{o}ng-Tr\`{i}nh L\^{e} }

\address{Department of Mathematics, Quy Nhon University\\
170 An Duong Vuong Street,  Quy Nhon City, Binh Dinh Province, Vietnam}

\email{lecongtrinh@qnu.edu.vn}
 \author{Trung-Hi\^eu Th\'{a}i }

\address{Department of Mathematics, Quy Nhon University\\
170 An Duong Vuong Street, Quy Nhon City, Binh Dinh Province, Vietnam}

\email{december2112@gmail.com}


\begin{abstract}
Two natural extensions of Jensen's functional equation on the real line are the equations $f(xy)+f(xy^{-1}) = 2f(x)$ and $f(xy)+f(y^{-1}x) = 2f(x)$, where $f$ is a map from a multiplicative group $G$ into an abelian additive group $H$. In a series of papers \cite{Ng1}, \cite{Ng2}, \cite{Ng3},  C. T. Ng has solved these functional equations for the case where $G$ is a free group and the linear group $GL_n(R)$, $R=\z,\r$, a quadratically closed field or a finite field. He has also mentioned, without detailed proof, in the above papers and in \cite{Ng4} that when $G$ is the symmetric group $S_n$ the group of all solutions of  these functional equations coincides with the group of all homomorphisms  from $(S_n,\cdot)$ to $(H,+)$. The aim of this paper is to give an elementary and direct proof of this fact.

\end{abstract}
\maketitle


\section{Introduction}
On the real line Jensen's functional equation can be written in the following form
$$f(x+y) + f(x-y) = 2f(x), \vm x, y \in \r. $$
Let $(G,\cdot)$ be a group with the neutral element $e$, $(H,+)$ an abelian group with the zero element $0$, $f: G \mtn H$ any map. The following natural extensions of Jensen's functional equation were considered by C.T. Ng (\cite{Ng1}, \cite{Ng2}, \cite{Ng3}):

\equ \label{pt1.1}
f(xy) + f(xy^{-1}) = 2f(x), \mbox{ for all } x,y \in G;
\eequ
\equ \label{pt1.2}
f(xy) + f(y^{-1}x) = 2f(x), \mbox{ for all } x,y \in G.
\eequ
By considering $g(x):=f(x) - f(e)$ for all $x\in G$ we may assume that
\equ \label{pt1.3}
f(e) = 0.
\eequ
Denote by $S_1(G,H)$ resp. $S_2(G,H)$ the set of all solutions of the functional equation (\ref{pt1.1}) resp. (\ref{pt1.2}) with the normalized condition (\ref{pt1.3}). Denote by $\Hom(G,H)$ the set of all homomorphisms from $G$ to $H$. These sets  are abelian additive groups and it is clear that
\equ\label{pt1.4}
\Hom(G,H)\leq S_1(G,H)
\eequ
and
\equ\label{pt1.5}
\Hom(G,H) \leq S_2(G,H).
\eequ
The equalities in (\ref{pt1.4}) and (\ref{pt1.5}) occur just in some special cases, for example, when $G$ is abelian and $H$ has no element of order $2$; $G$ a free group; $G$ the linear group $GL_n(R), n\geq 2$, $R=\z,\r$, a quadratically closed field or a finite field. For the case where $G$ is the symmetric group $S_n, n\geq 1$, C.T. Ng mentioned in \cite{Ng3} and \cite{Ng4} that the above equalities  also hold, however the author has not given a direct proof for this case. In this paper we give an elementary and direct  proof for the equalities in (\ref{pt1.4}) and (\ref{pt1.5}) when $G=S_n, n\geq 1$.

\section{The functional equation $f(xy) + f(xy^{-1}) = 2f(x)$}
In this section we consider the functional equation (\ref{pt1.1}) with  the normalized condition (\ref{pt1.3}). Denote by $S_n, n\geq 1$, the symmetric group on $n$ elements. The following theorem is the  main result in this section.
\thr \label{dl2.1} $S_1(S_n,H) = \Hom(S_n,H)$.
\ethr
To prove Theorem \ref{dl2.1} we need the following formulae.
\pro[{\cite[Theorem 2]{Ng1}}] \label{md2.1} Let $(G,\cdot)$ and $(H,+)$ be groups. Then for each $f\in S_1(G,H)$ and for all $n\in \z, x,y,z \in G$ we have
\equ\label{pt2.2}
f(xyz)+f(xzy)  = 2f(xy) + 2f(xz) -2f(x);
\eequ
\equ\label{pt2.3}
f(xyz)+f(yxz) = 2f(xz) + 2f(yz) - 2f(z);
\eequ
\equ\label{pt2.4}
2f(xyz) = 2f(xy) + 2f(xz) + 2f(yz) - 2f(x) - 2f(y) - 2f(z);
\eequ
\equ\label{pt2.5}
f(xyz) - f(xzy) = 2f(yz) - 2f(y) - 2f(z);
\eequ
\equ \label{pt2.1}
f(xy^nz) = nf(xyz) - (n-1)f(xz).
\eequ
In particular, we have
\equ \label{pt2.8}
f(x^n) = nf(x), \mbox{ for all } x \in G.
\eequ
\epro

\lm \label{lm2.1} Let $\sigma = (a~ b)$ and $\tau = (b ~c)$ be transpositions in $S_n$ ($a \not = c$). Then for each $f\in S_1(S_n,H)$ we have
$$f(\sigma\tau) = f(\sigma) + f(\tau)=0.$$
\elm
\pf[Proof] Denote by  $\delta$ the transposition $(a ~c)$. It is easy to verify that
\equ\label{pt2.9}
(\sigma \tau)^4 = \sigma \tau;
\eequ
\equ\label{pt2.10}
\begin{array}{c}
\sigma \delta \tau = \tau \delta \sigma = \delta;\\
\tau \sigma \delta = \delta \sigma \tau = \sigma;\\
\delta \tau \sigma = \sigma \tau \delta = \tau.
\end{array}
\eequ
Substituting $x$ by $\delta$, $y$ by $\sigma \tau$ in (\ref{pt1.1}), noting that $(\sigma\tau)^{-1}=\tau\sigma$, we have
$$ f(\delta \sigma \tau) + f(\delta \tau \sigma) = 2 f(\delta).$$
It follows from (\ref{pt2.10}) and (\ref{pt2.8})  that
$$ f(\sigma) + f(\tau)  =   f(\delta^2) = f(e) = 0.$$
On the other hand, substituting $x$ by $\sigma \tau$, $y$ by $\delta$ in (\ref{pt1.1}), noting that $\delta^{-1} = \delta$, we have
$$ f(\sigma \tau\delta) + f(\sigma\tau \delta) = 2 f(\sigma \tau).$$
It follows from (\ref{pt2.9}), (\ref{pt2.8}) and (\ref{pt2.10}) that
$$f(\sigma \tau) = f\big((\sigma\tau)^4\big) = 4f(\sigma \tau) = 2[2f(\sigma\tau)] = 2[2f(\tau)] = 2f(\tau^2) = 0.   $$
Therefore, $f(\sigma \tau) = 0 = f(\sigma) + f(\tau)$.
\epf

\lm\label{lm2.2} Let $\sigma = (a~ b)$ and $\tau = (c ~d)$ be transpositions in $S_n$ ($a, b, c $ and $d$ distinguished). Then for each $f\in S_1(S_n,H)$ we have
$$f(\sigma\tau) = f(\sigma) + f(\tau)=0.$$
\elm
\pf[Proof]
It follows from (\ref{pt2.10}) that $\sigma \tau = (a~b)[(a~c)(c~d)(a~d)]$. Substituting $x$ by $(a~b)(a~c)$, $y$ by $(c~d)(a~d)$ in (\ref{pt1.1}), noting that $\big((c~d)(a~d)\big)^{-1} = (a~d)(c~d)$ and $(a~b)(a~c)(a~d)(c~d) = (a~d)(d~b)$, we have
$$ f(\sigma\tau) + f\big((a~d)(d~b)\big) = 2f\big((a~b)(a~c)\big).$$
It follows from Lemma \ref{lm2.1} that
$$f\big((a~d)(d~b)\big) = f\big((a~b)(a~c)\big)  = 0.$$
Therefore $f(\sigma\tau) = 0$. \\
On the other hand, it  follows also  from Lemma \ref{lm2.1} that
$$ f\big((a~b)\big) +  f\big((b~c)\big) = 0; $$
$$  f\big((b~c)\big) +  f\big((c~d)\big) = 0.$$
Taking the summation of these two equations, noting that $2 f\big((b~c)\big) =  f\big((b~c)^2\big) = f(e) = 0$, we obtain
$$ f(\sigma) + f(\tau) = 0.$$
Therefore, $f(\sigma\tau) = 0 = f(\sigma) + f(\tau)$.
\epf

\lm \label{lm2.3} The product of two arbitrary  transpositions in $S_n$ has always a square root.
\elm
\pf[Proof] Let $\sigma$ and $\tau$ be transpositions in $S_n$. If $\sigma = (a~b) $ and $\tau = (b~c)$, it follows from (\ref{pt2.9}) that
$$ \sigma \tau = (\sigma\tau)^4 = [(\sigma\tau)^2]^2.$$
On the other hand, if $\sigma = (a~b)$ and $\tau = (c~d)$ then it is easy to verify that
$$\sigma \tau = (a~b)(c~d) = [(a~c)(c~b)(b~d)]^2.$$
\epf

\lm \label{lm2.4} For each $f\in S_1(S_n,H)$, $2f(x) = 0$ for all $x\in S_n$.
\elm
\pf[Proof]
The statement is trivial for $S_1=\{e\}$. For $S_2$, $x^2=e$ for all $x\in S_2$. Therefore it follows from (\ref{pt2.8}) that
$$2f(x) = f(x^2) = f(e)=0, \vm x \in S_2.$$
Now we consider the case where $n\geq 3$.
Since each permutation in $S_n$  can be written as a product of transpositions (see, for example, \cite[Corollary 1, p. 293]{Su}), for each $x\in S_n$,
$$ x = \sigma_1 \sigma_2\cdots \sigma_r, ~ \mbox{each } \sigma_i \mbox{ is a transposition in  } S_n.$$
We prove the lemma by induction on $r$. \\
For $r=1$, it follows from (\ref{pt2.8}) that $2f(x)=f(x^2) = f(e)=0$.  Assume that the lemma holds for $r \geq 1$.  We show that the lemma also holds for $r+1$.  In fact, it follows from (\ref{pt2.4}) that
\begin{align}
2f(x) &= 2  f\big((\sigma_1\cdots \sigma_{r-1})\cdot \sigma_r \cdot \sigma_{r+1}\big)  \nonumber\\
& = 2  f\big((\sigma_1\cdots \sigma_{r-1})\cdot \sigma_r\big) + 2  f\big((\sigma_1\cdots \sigma_{r-1})\cdot \sigma_{r+1}\big)+\nonumber\\
& + 2f(\sigma_r\sigma_{r+1}) - 2f(\sigma_1\cdots \sigma_{r-1}) - 2 f(\sigma_r) - 2 f(\sigma_{r+1}).\nonumber
\end{align}
The conclusion for $r+1$ follows from  Lemma \ref{lm2.1}, Lemma \ref{lm2.2} and the induction hypothesis.
\epf

\coro \label{coro2.1} For each $f\in S_1(S_n,H)$ and for every $x,y,z \in S_n$, we have
$$f(xyz) = f(xzy) = f(yxz). $$
 Therefore we may re-arrange the order of the transpositions in each permutation  $x \in  S_n$  in such a way that the value $f(x)$ does not change.
\ecoro
\pf[Proof] It follows from (\ref{pt2.2}), (\ref{pt2.3}), (\ref{pt2.5}) and Lemma \ref{lm2.4} that
$$ f(xyz) + f(yxz) = f(xyz) + f(xzy) = f(xyz) - f(xzy) = 0.$$
This proves the corollary.
\epf

\lm \label{lm2.5} For each $f\in S_1(S_n,H)$ and for each $x\in S_n$, the following holds:
\it
\item[(i)] If the permutation $x$ is even then $f(x) = 0$.
\item[(ii)] If $x$ is odd then $f(x) = -f(\sigma_r)$, where $\sigma_r$ is the last transposition in a decomposition of $x$.
\hit
\elm
\pf[Proof]
Let  $x\in S_n$. Then $x$ can be written as  $x=\sigma_1\cdots \sigma_r$, where each $\sigma_i$  is a transposition.  If the number of transpositions $r$ is even, namely $r=2k$, it follows from Lemma \ref{lm2.3} that there exist permutations $\tau_1,\cdots,\tau_{k}$ such that
$$x=\tau_1^2\cdots \tau_k^2.$$
Then it follows from Lemma \ref{lm2.4} that
$$ f(x) = f\big((\tau_1\cdots\tau_k)^2\big)  = 2 f(\tau_1\cdots \tau_k) = 0.$$
If $r$ is odd, namely $r=2k+1$, we can write  $x=\tau_1^2\cdots \tau_k^2 \sigma_{r}$. Then it follows from (\ref{pt2.1}) and Lemma  \ref{lm2.4} that
$$f(x) = f\big((\tau_1\cdots \tau_k)^2\sigma_{r}\big) = 2f\big((\tau_1\cdots \tau_k)\sigma_{r}\big) - f(\sigma_{r}) =  - f(\sigma_{r}).  $$
\epf

\pf[\textbf{Proof of Theorem \ref{dl2.1}}] Let $f$ be an arbitrary element of $S_1(S_n,H)$,  $x$ and $y$  any permutations in $S_n$. Then these permutations can be written as $x = \sigma_1\cdots\sigma_r$ and $y=\tau_1\cdots \tau_s$, where each $\sigma_i$ and $\tau_j$ are transpositions in $S_n$. We have the following cases:\\
\textbf{The first case:} both $r$ and $s$ are even. Then the product $xy$ is an even permutation, hence it follows from Lemma \ref{lm2.5} that $f(x)=f(y)=f(xy) = 0$. Thus $f(xy)=f(x)+f(y)$.\\
\textbf{The second case:} $r$ is even,   $s$ is odd. Then the product $xy$ is odd, hence it follows from Lemma \ref{lm2.5} that
$f(x)=0, f(y)=-f(\tau_s) $ and
$$ f(xy) = f\big((\sigma_1\cdots \sigma_r)(\tau_1\cdots \tau_{s-1})\tau_s\big) = -f(\tau_s) = f(x)+f(y).$$
\textbf{The third case:} $r$ is odd, $s$ is even. Similarly to the second case, we have
$$ f(xy) = -f(\sigma_r) = f(x) + f(y).$$
\textbf{The last case:} both $r$ and $s$ are odd. Then  $f(x)=-f(\sigma_r), f(y)=-f(\tau_s)$, while  $f(xy)=0$ since the product $xy$ is even. It follows from the proof of Lemma \ref{lm2.1} and Lemma \ref{lm2.2} that $f(\sigma_r)+f(\tau_s) = 0$. Therefore in this case we have also $f(xy) = 0 = f(x) + f(y)$. \\
Thus in any cases we always have $f(xy) = f(x) + f(y)$, i.e., $f\in \Hom(S_n,H)$. This proves the theorem.
\epf

\section{The functional equation $f(xy) + f(y^{-1}x) = 2f(x)$}
In this section we consider the functional equation (\ref{pt1.2}) with  the normalized condition (\ref{pt1.3}) and show that the equality in (\ref{pt1.5}) occurs  for the group $G=S_n$. The proof for this equality follows step-by-step the one given for the equation (\ref{pt1.1}) in Section 2.
\pro[{\cite[Theorem 2.1]{Ng3}}] \label{md3.1}
Let $(G,\cdot)$ and $(H,+)$ be groups. Then for each $f\in S_2(G,H)$ and for all $n\in \z, x,y,z \in G$ we have
\equ \label{pt3.1}
f(x^n) = nf(x);
\eequ
\equ\label{pt3.2}
f(xyz)+f(xzy)  = 2f(xy) + 2f(xz) -2f(x);
\eequ
\equ\label{pt3.3}
f(xyz)+f(yxz) = 2f(xz) + 2f(yz) - 2f(z);
\eequ
\equ\label{pt3.4}
f(xyz) - f(xzy) = 2f(yz) - 2f(y) - 2f(z);
\eequ
\equ\label{pt3.5}
2f(xyz) = 2f(xy) + 2f(xz) + 2f(yz) - 2f(x) - 2f(y) - 2f(z);
\eequ
\equ \label{pt3.6}
f(xy^2z) = f(xz) + 2f(y).
\eequ
\epro

\lm \label{lm3.1}
Let $\sigma = (a~ b)$ and $\tau = (b ~c)$ be transpositions in $S_n$ ($a \not = c$). Then for each $f\in S_2(S_n,H)$ we have
$$f(\sigma\tau) = f(\sigma) + f(\tau)=0.$$
\elm
\pf[Proof] Denote $\delta = (a~c)$. Substituting $x$ by $\sigma\tau$ and $y$ by $\delta$ in (\ref{pt1.2}), noting that $\delta^{-1}=\delta$ and using (\ref{pt2.10}) we have
\equ \label{pt3.7}
 f(\tau) + f(\sigma) = 2 f(\sigma\tau).
\eequ
It follows from (\ref{pt3.7}),  (\ref{pt2.9}) and (\ref{pt3.1}) that
\begin{align}
 f(\sigma\tau) &= f\big((\sigma\tau)^4\big) = 4f(\sigma\tau)= \nonumber \\
 &= 2f(\sigma)+2f(\tau) = f(\sigma^2)+f(\tau^2) = f(e)+f(e)=0.\nonumber
\end{align}
Substituting $f(\sigma\tau)=0$ into the equation (\ref{pt3.7}) we have
$$f(\sigma) + f(\tau) = 2f(\sigma\tau) = 0 = f(\sigma\tau).$$
\epf

\lm\label{lm3.2} Let $\sigma = (a~ b)$ and $\tau = (c ~d)$ be transpositions in $S_n$ ($a, b, c $ and $d$ distinguished). Then for each $f\in S_2(S_n,H)$ we have
$$f(\sigma\tau) = f(\sigma) + f(\tau)=0.$$
\elm
\pf[Proof] Substituting $x  = (a~b)(a~c)$ and $y=(c~d)(a~d)$ in (\ref{pt1.2}), noting that
$$xy=(a~b)(a~c)(c~d)(a~d) =  (a~b)[(a~c)(c~d)(a~d)] \overset{(\ref{pt2.10})}{==} (a~b)(c~d) = \sigma \tau,$$
$$ \big((c~d)(a~d)\big)^{-1} = (a~d)(c~d),$$
and $y^{-1}x=(a~d)(c~d) (a~b)(a~c) = (b~d)(d~c)$, we have
$$f(\sigma\tau) + f\big((b~d)(d~c)\big) = 2f\big((a~b)(a~c)\big). $$
Hence it follows from Lemma \ref{lm3.1} that $f(\sigma\tau) = 0$. \\
By a similar argument to the proof of the second part of Lemma \ref{lm2.2} we have also $f(\sigma)+f(\tau)=0=f(\sigma\tau)$.
\epf

\lm \label{lm3.3}
For each $f\in S_2(S_n,H)$, $2f(x) = 0$ for all $x\in S_n$.
\elm
\pf[Proof]
Since we have also the equation (\ref{pt3.5}) which is the same as the equation (\ref{pt2.4}), the proof for this lemma is the same as the one given in Lemma \ref{lm2.4}.
\epf

\coro \label{coro3.1} For each $f\in S_2(S_n,H)$ and for every $x,y,z \in S_n$, we have
$$f(xyz) = f(xzy) = f(yxz). $$
 Therefore we may re-arrange the order of the transpositions in each permutation  $x \in  S_n$  in such a way that the value $f(x)$ does not change.
\ecoro
\pf[Proof] The corollary  follows from (\ref{pt3.2}), (\ref{pt3.3}), (\ref{pt3.4}) and Lemma \ref{lm3.3}.
\epf

\lm \label{lm3.4} For each $f\in S_2(S_n,H)$ and for each $x\in S_n$, the following holds:
\it
\item[(i)] If the permutation $x$ is even then $f(x) = 0$.
\item[(ii)] If $x$ is odd then $f(x) = f(\sigma_r)$, where $\sigma_r$ is the last transposition in a decomposition of $x$.
\hit
\elm
\pf[Proof] The proof for (i) is the same as  the one given in Lemma \ref{lm2.5} (i), using Lemma \ref{lm3.3}.
The proof for (ii) is also similar to the one given in Lemma \ref{lm2.5} (ii),  using the equation (\ref{pt3.6}) instead of the equation (\ref{pt2.1}).
\epf
Therefore, by a similar argument to the proof of Theorem \ref{dl2.1}, using Lemma \ref{lm3.4} instead of Lemma \ref{lm2.5}, we obtain the main theorem of this section.

\thr \label{dl3.1}
$S_2(S_n,H) = \Hom(S_n,H)$.
\ethr
\textbf{Acknowledgement:} {The author is indebted to  the referees for their useful comments and suggestions. }

\end{document}